\documentclass[letterpaper, 10 pt, conference]{ieeeconf}  

\IEEEoverridecommandlockouts                              
\usepackage{amsmath} 
\usepackage{amssymb}  
\usepackage{color,xcolor}
\title{\LARGE \bf
Connection between MP and DPP for Stochastic Recursive Optimal Control Problems: Viscosity Solution Framework in Local Case
}

\author{Tianyang Nie, Jingtao Shi, {\it Member, IEEE}, Zhen Wu
\thanks{This work was supported by the National Natural Science Foundations of China (11301011, 11201264, 11571205, 61573217), the 111 project (B12023), the Chang Jiang Scholar Program of Chinese Education Ministry, the Natural Science Foundations of Shandong Province of China (JQ201401, ZR2015JL003) and the Fundamental Research Fund of Shandong University (2015HW023).}
\thanks{Corresponding author: Jingtao Shi. \ T.~Y. Nie, J.~T. Shi and Z. Wu are all with School of Mathematics, Shandong University, Jinan 250100, P.~R. China.
        {\tt\small nietianyang@sdu.edu.cn, shijingtao@sdu.edu.cn, wuzhen@sdu.edu.cn}}%
}

\begin{document}

\maketitle
\thispagestyle{empty}
\pagestyle{empty}

\begin{abstract}

This paper deals with a nonsmooth version of the connection between the maximum principle and dynamic programming principle, for the stochastic recursive control problem when the control domain is convex. By employing the notions of sub- and super-jets, the set inclusions are derived among the value function and the adjoint processes. The general case for non-convex control domain is open.

\end{abstract}

\section{INTRODUCTION}

\noindent

There are usually two ways to study optimal control problems: Pontryagin's maximum principle (MP) and Bellman's dynamic programming principle (DPP), involving an adjoint variable $\psi$ and the value function $V$, respectively. The classical result by Fleming and Rishel \cite{FR75} on the connection between the MP and DPP is known as $\psi(t)=-V_x(t,\bar{x}(t))$, where $\bar{x}(\cdot)$ is the optimal state. Since the value function $V$ is not always smooth, some non-smooth versions of the classical result are researched by non-smooth analysis and generalized derivatives. Within the framework of viscosity solution, Zhou \cite{Zhou90} showed that
\begin{equation}\label{Zhou90}
\begin{aligned}
D_x^{1,-}V(t,\bar{x}(t))\subset\{-\psi(t)\}\subset D_x^{1,+}V(t,\bar{x}(t)),
\end{aligned}
\end{equation}
where $D_x^{1,-}V(t,\bar{x}(t))$ and $D_x^{1,+}V(t,\bar{x}(t))$ denote the first-order sub- and super-jets of $V$ at $(t,\bar{x}(t))$, respectively.

For stochastic optimal control problems, the classical result on the connection between the MP and DPP is proved by Bensoussan \cite{Ben82}, which is known as $p(t)=-V_x(t,\bar{x}(t)),\ q(t)=-V_{xx}(t,\bar{x}(t))\sigma(t,\bar{x}(t),\bar{u}(t))$ involving an adjoint process pair $(p,q)$, where $\bar{u}(\cdot)$ is the optimal control and $\sigma$ is the diffusion coefficient. Within the framework of viscosity solution, Yong and Zhou \cite{YZ99} showed that
\begin{equation}\label{YZ99}
\begin{aligned}
&\{-p(t)\}\times[-P(t),\infty)\subset D_x^{2,+}V(t,\bar{x}(t)),\\
&D_x^{2,-}V(t,\bar{x}(t))\subset \{-p(t)\}\times(-\infty,-P(t)],
\end{aligned}
\end{equation}
where $D_x^{2,-}V(t,\bar{x}(t))$ and $D_x^{2,+}V(t,\bar{x}(t))$ denote the second-order sub- and super-jets of $V$ at $(t,\bar{x}(t))$, and $p,P$ are the first- and second-order adjoint processes, respectively.

In this paper, we consider one kind of stochastic recursive optimal control problem, where the cost functional is described by the solution to a {\it backward stochastic differential equation} (BSDE) of the following form
\begin{equation*}
\left\{
\begin{aligned}
-dy(t)&=f(t,y(t),z(t))dt-z(t)dW(t),\ t\in[0,T],\\
  y(T)&=\xi,
\end{aligned}
\right.
\end{equation*}
where the terminal condition (rather than the initial condition) $\xi$ is given in advance.
Linear BSDE was introduced by Bismut \cite{Bis78}, to represent the adjoint equation when applying the MP to solve stochastic optimal control problems. The nonlinear BSDE was introduced by Pardoux and Peng \cite{PP90}. Independently, Duffie and Epstein \cite{DE92} introduced BSDE from economic background, and they presented a stochastic differential formulation of recursive utility which is an extension of the standard additive utility with the instantaneous utility depending not only on the instantaneous consumption rate but also on the future utility. Stochastic recursive optimal control problems have found important applications in mathematical economics, mathematical finance and engineering (see El Karoui, Peng and Quenez \cite{EPQ97,EPQ01}, Wang and Wu \cite{WW09}, Cvitanic and Zhang \cite{CZ13} and the references therein).

For stochastic recursive optimal control problems, Peng \cite{Peng93} first obtained a local maximum principle when the control domain is convex. And Xu \cite{Xu95} studied the non-convex control domain case, but with the assumption that the diffusion coefficient does not depend on the control variable. Wu \cite{Wu13} established a general maximum principle by Ekeland variational principle, where the control domain is non-convex and the diffusion coefficient contains the control variable. Peng \cite{Peng92} (also see Peng \cite{Peng97}) first obtained the generalized dynamic programming principle and introduced a generalized {\it Hamilton-Jacobi-Bellman} (HJB) equation which is a second-order parabolic {\it partial differential equation} (PDE). The value function is proved to be the viscosity solution to the generalized HJB equation.

The connection between MP and DPP for stochastic recursive optimal control problems was first studied by Shi \cite{Shi10} (see also Shi and Yu \cite{SY13}) in its local form, when the control domain is convex and the value function is assumed to be smooth enough. The main result is
{\small\begin{equation}\label{Shi10}
\left\{
\begin{aligned}
p(t)&=V_x(t,\bar{x}(t))^\top q(t),\\
k(t)&=\big[V_{xx}(t,\bar{x}(t))\sigma(t,\bar{x}(t),\bar{u}(t))+V_x(t,\bar{x}(t))f_z\big(t,\bar{x}(t),\\
    &\quad -V(t,\bar{x}(t)),-V_x(t,\bar{x}(t))\sigma(t,\bar{x}(t),\bar{u}(t)),\bar{u}(t)\big)\big]q(t),
\end{aligned}
\right.
\end{equation}}
involving an adjoint process triple $(p,q,k)$, where $f$ is the generator of the controlled BSDE which is coupled with the controlled SDE. Applications to the recursive utility portfolio optimization problem in the financial market are discussed.

However, this classical result is highly unsatisfactory because the smoothness assumption on the value function $V$ is illusory and it is not true even in the very simple case: see Example 3.1 of this paper. In the current work, we extend the above classical result by getting rid of the illusory assumption that the value function is differentiable. Our main contribution is to show the connection between the adjoint processes $p(\cdot),q(\cdot)$ in the maximum principle and the first-order sub- and super-jets $D_x^{1,-}V(t,\bar{x}(t)),D_x^{1,+}V(t,\bar{x}(t))$.

The rest of this paper is organized as follows. In Section 2, we state our problem and give some preliminary results about the MP and the DPP. Section 3 exhibits the main result of this paper, namely, the connection between the value function and the adjoint processes within the framework of viscosity solution. Finally, in Section 4 we give the concluding remarks.

\section{Problem Statement and Preliminaries}

\noindent

Let $T>0$ be finite and $\mathbf{U}\subset\mathbf{R}^k$ be nonempty and convex. Given $t\in[0,T)$, we denote $\mathcal{U}^w[t,T]$ the set of all 5-tuples $(\Omega,\mathcal{F},\mathbf{P},W(\cdot);u(\cdot))$ satisfying the following:

(i) $(\Omega,\mathcal{F},\mathbf{P})$ is a complete probability space;

(ii) $\{W(s)\}_{s\geq t}$ is a $d$-dimensional standard Brownian motion defined on $(\Omega,\mathcal{F},\mathbf{P})$ over $[t,T]$ (with $W(t)=0$ almost surely), and
$\mathcal{F}_s^t=\sigma\{W(r);t\leq r\leq s\}$ augmented by all the $\mathbf{P}$-null sets in $\mathcal{F}$;

(iii) $u:[t,T]\times\Omega\rightarrow\mathbf{U}$ is an $\{\mathcal{F}_s^t\}_{s\geq t}$-adapted process on $(\Omega,\mathcal{F},\mathbf{P})$.

We write $(\Omega,\mathcal{F},\mathbf{P},W(\cdot);u(\cdot))\in\mathcal{U}^w[t,T]$, but occasionally we will write only $u(\cdot)\in\mathcal{U}^w[t,T]$ if no ambiguity exists. For any $(t,x)\in[0,T)\times\mathbf{R}^n$, consider the state $X^{t,x;u}(\cdot)\in{\mathbf{R}}^n$ given by the following controlled SDE:
{\small\begin{equation}\label{controlled SDE}
\left\{
\begin{aligned}
 dX^{t,x;u}(s)&=b(s,X^{t,x;u}(s),u(s))ds\\
              &\quad+\sigma(s,X^{t,x;u}(s),u(s))dW(s),\ s\in[t,T],\\
  X^{t,x;u}(t)&=x.
\end{aligned}
\right.
\end{equation}}
Here $b:[0,T]\times\mathbf{R}^n\times\mathbf{U}\rightarrow\mathbf{R}^n,
\sigma:[0,T]\times\mathbf{R}^n\times\mathbf{U}\rightarrow\mathbf{R}^{n\times d}$ are given functions. We assume that

\vspace{1mm}

\noindent({\bf H1})\quad $b,\sigma$ are uniformly continuous in $(s,x,u)$, and there exists a constant $C>0$ such that for all
$s\in[0,T],x,\hat{x}\in\mathbf{R}^n,u\in\mathbf{U}$,
{\small\begin{equation*}
\left\{
\begin{aligned}
 &|b(s,x,u)-b(s,\hat{x},u)|+|\sigma(t,x,u)-\sigma(s,\hat{x},u)|\leq C|x-\hat{x}|,\\
 &|b(s,x,u)|+|\sigma(s,x,u)|\leq C(1+|x|).
\end{aligned}
\right.
\end{equation*}}
For any $u(\cdot)\in\mathcal{U}^w[t,T]$, under ({\bf H1}), SDE (\ref{controlled SDE}) has a unique solution $X^{t,x;u}(\cdot)$ by the classical SDE theory (see \cite{N97,YZ99}). We refer to such
$u(\cdot)\in\mathcal{U}^w[t,T]$ as an admissible control and $(X^{t,x;u}(\cdot),u(\cdot))$ as an admissible pair.

Next, we introduce the following controlled BSDE coupled with (\ref{controlled SDE}):
{\small\begin{equation}\label{controlled BSDE}
\left\{
\begin{aligned}
-dY^{t,x;u}(s)&=f(s,X^{t,x;u}(s),Y^{t,x;u}(s),Z^{t,x;u}(s),u(s))ds\\
              &\quad-Z^{t,x;u}(s)dW(s),\ s\in[t,T],\\
  Y^{t,x;u}(T)&=\phi(X^{t,x;u}(T)).
\end{aligned}
\right.
\end{equation}}
Here $f:[0,T]\times\mathbf{R}^n\times\mathbf{R}\times\mathbf{R}^d\times\mathbf{U}\rightarrow\mathbf{R},\Phi:\mathbf{R}^n\rightarrow\mathbf{R}$ are given functions. We assume that

\vspace{1mm}

\noindent({\bf H2})\quad $f,\phi$ are uniformly continuous in $(s,x,y,z,u)$ and there exists a constant $C>0$
such that for all $s\in[0,T],x,\hat{x}\in\mathbf{R}^n,y,\hat{y}\in\mathbf{R},z,\hat{z}\in\mathbf{R}^d,u\in\mathbf{U}$,
\begin{equation*}
\left\{
\begin{aligned}
 &|f(s,x,y,z,u)-f(s,\hat{x},\hat{y},\hat{z},u)|\\
 &\ \leq C(|x-\hat{x}|+|y-\hat{y}|+|z-\hat{z}|),\\
 &|f(s,x,0,0,u)|+|\phi(x)|\leq C(1+|x|),\\
 &|\phi(x)-\phi(\hat{x})|\leq C|x-\hat{x}|.
\end{aligned}
\right.
\end{equation*}
Then for any $u(\cdot)\in\mathcal{U}^w[t,T]$ and the given unique solution $X^{t,x;u}(\cdot)$ to (\ref{controlled SDE}), under ({\bf H2}), BSDE (\ref{controlled BSDE}) admits a unique solution $(Y^{t,x;u}(\cdot),Z^{t,x;u}(\cdot))$ by the classical BSDE theory (see Pardoux and Peng \cite{PP90} or Peng \cite{Peng97}).

Given $u(\cdot)\in\mathcal{U}^w[t,T]$, we introduce the cost functional
\begin{equation}\label{cost functional}
J(t,x;u(\cdot)):=-Y^{t,x;u}(s)|_{s=t},\quad(t,x)\in[0,T]\times\mathbf{R}^n.
\end{equation}

Our recursive stochastic optimal control problem is the following.

\vspace{1mm}

\noindent{\bf Problem (RSOCP).}\quad For given $(t,x)\in[0,T)\times\mathbf{R}^n$, to minimize (\ref{cost
functional}) subject to (\ref{controlled SDE})$\sim$(\ref{controlled BSDE}) over $\mathcal{U}^w[t,T]$.

We define the value function
{\small\begin{equation}\label{value function}
\left\{
\begin{aligned}
&V(t,x):=\inf\limits_{u(\cdot)\in\mathcal{U}^w[t,T]}J(t,x;u(\cdot)),\ (t,x)\in[0,T]\times\mathbf{R}^n,\\
&V(T,x)=-\phi(x),\quad x\in\mathbf{R}^n.
\end{aligned}
\right.
\end{equation}
Any $\bar{u}(\cdot)\in\mathcal{U}^w[t,T]$ that achieves the above infimum is called an optimal control, and the corresponding solution triple $(\bar{X}^{t,x;\bar{u}}(\cdot),\bar{Y}^{t,x;\bar{u}}(\cdot),\bar{Z}^{t,x;\bar{u}}(\cdot))$ is called an optimal state. We refer to $(\bar{X}^{t,x;\bar{u}}(\cdot),\bar{Y}^{t,x;\bar{u}}(\cdot),\bar{Z}^{t,x;\bar{u}}(\cdot),\bar{u}(\cdot))$ as an optimal quadruple.

\vspace{1mm}

\noindent{\bf Remark 2.1}\quad Because $b,\sigma,f,g$ are all deterministic functions, then from Proposition 5.1 of Peng \cite{Peng97}, we know that under ({\bf H1}), ({\bf H2}), the above value function is a deterministic function. Thus our definition (\ref{value function}) is meaningful.

We introduce the following generalized HJB equation:
\begin{equation}\label{HJB equation}
\left\{
\begin{aligned}
 &-v_t(t,x)+\sup\limits_{u\in\mathbf{U}}G\big(t,x,-v(t,x),-v_x(t,x),\\
 &\quad-v_{xx}(t,x),u\big)=0,\ (t,x)\in[0,T)\times\mathbf{R}^n,\\
 &v(T,x)=-\phi(x),\quad \forall x\in\mathbf{R}^n,
\end{aligned}
\right.
\end{equation}
where the generalized Hamiltonian function $G:[0,T]\times\mathbf{R}^n\times\mathbf{R}\times\mathbf{R}^n\times\mathcal{S}^n\times\mathbf{U}\rightarrow\mathbf{R}$ is defined as
\begin{equation}\label{generalized Hamiltonian}
\begin{aligned}
 &G(t,x,r,p,A,u):=\frac{1}{2}\mbox{tr}\big\{\sigma(t,x,u)^\top A\sigma(t,x,u)\big\}\\
               &\quad+\langle p,b(t,x,u)\rangle+f(t,x,r,\sigma(t,x,u)^\top p,u).
\end{aligned}
\end{equation}

The following result belongs to Peng \cite{Peng97}.

\vspace{1mm}

\noindent{\bf Proposition 2.1}\quad\quad{\it Let {\bf (H1), (H2)} hold. Then for any $t\in[0,T]$ and $x,x'\in\mathbf{R}^n$, we have
\begin{equation}\label{value function: regularity}
\begin{aligned}
&\mbox{(i)}\quad |V(t,x)-V(t,x')|\leq C|x-x'|,\\
&\mbox{(ii)}\quad |V(t,x)|\leq C(1+|x|).
\end{aligned}
\end{equation} }

\vspace{1mm}

We introduce the definition of the viscosity solution for HJB equation (\ref{HJB equation}).

\vspace{1mm}

\noindent{\bf Definition 2.1}\quad{\it (i) A function $v\in C([0,T]\times\mathbb{R}^n)$ is called a viscosity subsolution to (\ref{HJB equation}) if
\begin{equation*}
v(T,x)\leq -\phi(x),\quad \forall x\in\mathbb{R}^n,
\end{equation*}
and for any $\varphi\in C^{1,2}([0,T]\times\mathbb{R}^n)$, whenever $v-\varphi$ attains a local maximum at $(t,x)\in[0,T)\times\mathbb{R}^n$, we have
{\small\begin{equation*}
\begin{aligned}
-\varphi_t(t,x)+\sup\limits_{u\in\mathbf{U}}G\big(t,x,-v(t,x),-\varphi_x(t,x),-\varphi_{xx}(t,x),u\big)\leq0.
\end{aligned}
\end{equation*}}
(ii) A function $v\in C([0,T]\times\mathbb{R}^n)$ is called a viscosity supersolution to (\ref{HJB equation}) if
\begin{equation*}
v(T,x)\geq -\phi(x),\quad \forall x\in\mathbb{R}^n,
\end{equation*}
and for any $\varphi\in C^{1,2}([0,T]\times\mathbb{R}^n)$, whenever $v-\varphi$ attains a local minimum at $(t,x)\in[0,T)\times\mathbb{R}^n$, we have
{\small\begin{equation*}
\begin{aligned}
-\varphi_t(t,x)+\sup\limits_{u\in\mathbf{U}}G\big(t,x,-v(t,x),-\varphi_x(t,x),-\varphi_{xx}(t,x),u\big)\geq0.
\end{aligned}
\end{equation*}}
(iii) A function $v\in C([0,T]\times\mathbb{R}^n)$ is called a viscosity solution to (\ref{HJB equation}) if it is both a viscosity subsolution and viscosity supersolution to (\ref{HJB equation}).}

\vspace{1mm}

The following result also belongs to Peng \cite{Peng97}.

\vspace{1mm}

\noindent{\bf Proposition 2.2}\quad\quad{\it Let {\bf (H1), (H2)} hold. Then $V(\cdot,\cdot)$ defined by (\ref{value function}) is the unique viscosity solution to (\ref{HJB equation}).}

To conveniently state the maximum principle, we regard the above (\ref{controlled SDE}), (\ref{controlled BSDE}) as a controlled {\it forward-backward stochastic differential equation} (FBSDE):
{\small\begin{equation}\label{controlled FBSDE}
\left\{
\begin{aligned}
 dX^{t,x;u}(s)&=b(s,X^{t,x;u}(s),u(s))ds\\
              &\quad+\sigma(s,X^{t,x;u}(s),u(s))dW(s),\\
-dY^{t,x;u}(s)&=f(s,X^{t,x;u}(s),Y^{t,x;u}(s),Z^{t,x;u}(s),u(s))ds\\
              &\quad-Z^{t,x;u}(s)dW(s),\ s\in[t,T],\\
  X^{t,x;u}(t)&=x,\ Y^{t,x;u}(T)=\phi(X^{t,x;u}(T)).
\end{aligned}
\right.
\end{equation}}

We need the following assumption.

\vspace{1mm}

\noindent({\bf H3})\quad $b,\sigma,\phi,f$ are continuously differentiable in $(x,y,z)$ and the partial derivatives are uniformly bounded.

Let $(\bar{X}^{t,x;\bar{u}}(\cdot),\bar{Y}^{t,x;\bar{u}}(\cdot),\bar{Z}^{t,x;\bar{u}}(\cdot),\bar{u}(\cdot))$ be an optimal quadruple. For all $s\in[0,T]$, we denote
\begin{equation*}
\begin{aligned}
\bar{b}(s):=&\ b(s,\bar{X}^{t,x;\bar{u}}(s),\bar{u}(s)),\ \bar{\sigma}(s):=b(s,\bar{X}^{t,x;\bar{u}}(s),\bar{u}(s)),\\
\bar{f}(s):=&\ f(s,\bar{X}^{t,x;\bar{u}}(s),\bar{Y}^{t,x;\bar{u}}(s),\bar{Z}^{t,x;\bar{u}}(s),\bar{u}(s)),
\end{aligned}
\end{equation*}
and similar notations are used for all their derivatives.

We introduce the adjoint equation:
\begin{equation}\label{adjoint equation}
\left\{
\begin{aligned}
-dp(s)&=\big[\bar{b}_x(s)^\top p(s)-\bar{f}_x(s)^\top q(s)+\bar{\sigma}_x(s)k(s)\big]ds\\
      &\quad-k(s)dW(s),\\
 dq(s)&=\bar{f}_y(s)^\top q(s)ds+\bar{f}_z(s)^\top q(s)dW(s),\ s\in[t,T],\\
  p(T)&=-\phi_x(\bar{X}^{t,x;\bar{u}}(T))^\top q(T),\quad q(t)=1,
\end{aligned}
\right.
\end{equation}
and the Hamiltonian function
$H:[0,T]\times\mathbf{R}^n\times\mathbf{R}\times\mathbf{R}^d\times\mathbf{U}\times\mathbf{R}^n\times\mathbf{R}\times\mathbf{R}^{n\times d}\rightarrow\mathbf{R}$ is defined as
\begin{equation}\label{Hamiltonian}
\begin{aligned}
&H(t,x,y,z,u,p,q,k):=\langle p,b(t,x,u)\rangle\\
&-\langle q,f(t,x,y,z,u)\rangle+\mbox{tr}\big[\sigma(t,x,u)^\top k\big].
\end{aligned}
\end{equation}

Under {\bf (H1), (H2), (H3)}, (\ref{adjoint equation}) admits a unique solution $(p(\cdot),q(\cdot),k(\cdot))$, which is called the adjoint process triple.

The following result comes from Peng \cite{Peng93}.

\vspace{1mm}

\noindent{\bf Proposition 2.3}\quad{\it Let {\bf (H1), (H2), (H3)} hold
and $(t,x)\in[0,T)\times\mathbf{R}^n$ be fixed. Suppose that $\bar{u}(\cdot)$ is an optimal control for {\bf Problem (RSOCP)}, and
$(\bar{X}^{t,x;\bar{u}}(\cdot),\bar{Y}^{t,x;\bar{u}}(\cdot),\bar{Z}^{t,x;\bar{u}}(\cdot))$
is the corresponding optimal state. Let $(p(\cdot),q(\cdot),k(\cdot))$ be the adjoint process triple. Then
\begin{equation}\label{maximum condition}
\begin{aligned}
\big\langle
&H_u(s,\bar{X}^{t,x;\bar{u}}(s),\bar{Y}^{t,x;\bar{u}}(s),\bar{Z}^{t,x;\bar{u}}(s),\bar{u}(s),\\
&\quad p(s),q(s),k(s)),u-\bar{u}(s)\big\rangle\geq0,\ \forall u\in\mathbf{U},
\end{aligned}
\end{equation}
a.e. $s\in[t,T],\mathbf{P}\mbox{-}$a.s.}

\vspace{1mm}

\noindent{\bf Remark 2.2}\quad Notice that Proposition 2.3 is proved by Peng \cite{Peng93} in its strong formulation. However, as pointed out in Yong and Zhou \cite{YZ99}, since the DPP is involved, we need to deal with {\bf Problem (RSOCP)} in its weak formulation. Since only necessary conditions of optimality are considered here, an optimal quadruple (no matter whether in the strong or weak formulation) is given as a starting point, and all the results are valid for this given optimal quadruple on the probability space it attached to.

\section{Main Result}

\noindent

We first introduce the notion of the first-order super- and sub-jets. For $v\in C([0,T]\times\mathbb{R}^n)$, and $(t,\hat{x})\in[0,T]\times\mathbb{R}^n$, we define
\begin{equation}\label{first-order super- and sub-jets}
\left\{
\begin{aligned}
D_x^{1,+}v(t,\hat{x})&:=\Big\{p\in\mathbf{R}^n\big|v(t,x)\leq v(t,\hat{x})+\langle p,x-\hat{x}\rangle\\
                     &\qquad+o(|x-\hat{x}|),\mbox{ as }x\rightarrow\hat{x}\Big\},\\
D_x^{1,-}v(t,\hat{x})&:=\Big\{p\in\mathbf{R}^n\big|v(t,x)\geq v(t,\hat{x})+\langle p,x-\hat{x}\rangle\\
                     &\qquad+o(|x-\hat{x}|),\mbox{ as }x\rightarrow\hat{x}\Big\},
\end{aligned}
\right.
\end{equation}
\noindent{\bf Theorem 3.1}\quad{\it Let {\bf (H1), (H2), (H3)} hold
and $(t,x)\in[0,T)\times\mathbf{R}^n$ be fixed. Suppose that $\bar{u}(\cdot)$ is an optimal control for {\bf Problem (RSOCP)}, and
$(\bar{X}^{t,x;\bar{u}}(\cdot),\bar{Y}^{t,x;\bar{u}}(\cdot),\bar{Z}^{t,x;\bar{u}}(\cdot))$
is the corresponding optimal state. Let $(p(\cdot),q(\cdot),k(\cdot))$ be the adjoint process triple. Then
\begin{equation}\label{connection}
\begin{aligned}
&D_x^{1,-}V(s,\bar{X}^{t,x;\bar{u}}(s))\subset\{p(s)q^{-1}(s)\}\\
&\subset D_x^{1,+}V(s,\bar{X}^{t,x;\bar{u}}(s)),\ \forall s\in[t,T],\mathbf{P}\mbox{-}a.s.
\end{aligned}
\end{equation}
where $V(\cdot,\cdot)$ is the value function defined by (\ref{value function}).}

{\it Proof.}\quad Fix an $s\in[t,T]$. For any $x^1\in\mathbf{R}^n$, denote by $(X^{s,x^1;\bar{u}}(\cdot),Y^{s,x^1;\bar{u}}(\cdot),Z^{s,x^1;\bar{u}}(\cdot))$ the solution to the following FBSDE on $[s,T]$:
\begin{equation}\label{disturbed controlled FBSDE}
\left\{
\begin{aligned}
 X^{s,x^1;u}(r)&=x^1+\int_s^rb(\alpha,X^{s,x^1;u}(\alpha),u(\alpha))d\alpha\\
               &\quad+\int_s^r\sigma(\alpha,X^{s,x^1;u}(\alpha),u(\alpha))dW(\alpha),\\
 Y^{s,x^1;u}(r)&=\phi(X^{s,x^1;u}(T))+\int_r^Tf(\alpha,X^{s,x^1;u}(\alpha),\\
               &\qquad Y^{s,x^1;u}\alpha),Z^{s,x^1;u}(\alpha),u(\alpha))d\alpha\\
               &\quad-\int_r^TZ^{s,x^1;u}(\alpha)dW(\alpha),\ r\in[s,T].
\end{aligned}
\right.
\end{equation}
It is clear that (\ref{disturbed controlled FBSDE}) can be regarded as an FBSDE on $\big(\Omega,\mathcal{F},\{\mathcal{F}_r^t\}_{r\geq t},\mathbf{P}(\cdot|\mathcal{F}_s^t)(\omega)\big)$ for $\mathbf{P}\mbox{-}a.s. \omega$, where
$\mathbf{P}(\cdot|\mathcal{F}_s^t)(\omega)$ is the regular conditional probability given $\mathcal{F}_s^t$ defined on $(\Omega,\mathcal{F})$.

For any $s\leq r\leq T$, set
\begin{equation*}
\begin{aligned}
\hat{X}(r)&:=X^{s,x^1;\bar{u}}(r)-\bar{X}^{t,x;\bar{u}}(r),\\
\hat{Y}(r)&:=Y^{s,x^1;\bar{u}}(r)-\bar{Y}^{t,x;\bar{u}}(r),\\
\hat{Z}(r)&:=Z^{s,x^1;\bar{u}}(r)-\bar{Z}^{t,x;\bar{u}}(r).
\end{aligned}
\end{equation*}
Thus by a standard argument (see Theorem 6.3, Chapter 1, Yong and Zhou \cite{YZ99}), we have for any integer $k\geq1$,
\begin{equation}\label{estimate of SDE}
\mathbb{E}\Big[\sup\limits_{s\leq r\leq T}\big|\hat{X}(r)\big|^{2k}\Big|\mathcal{F}_s^t\Big]\leq C|x^1-\bar{X}^{t,x;\bar{u}}(s)|^{2k},\ \mathbf{P}\mbox{-}a.s.
\end{equation}
Moreover, the following estimates holds by Peng \cite{Peng97},
\begin{equation}\label{estimate of BSDE}
\begin{aligned}
\mathbb{E}\Big[\sup\limits_{s\leq r\leq T}\big|\hat{Y}(r)\big|^{2k}\Big|\mathcal{F}_s^t\Big]\leq C|x^1-\bar{X}^{t,x;\bar{u}}(s)|^{2k},\ \mathbf{P}\mbox{-}a.s.,\\
\mathbb{E}\Big[\Big(\int_s^T\big|\hat{Z}(r)\big|^2dr\Big)^k\Big|\mathcal{F}_s^t\Big]\leq C|x^1-\bar{X}^{t,x;\bar{u}}(s)|^{2k},\ \mathbf{P}\mbox{-}a.s.
\end{aligned}
\end{equation}
Now we write the equation for $\hat{X}(\cdot)$ as
\begin{equation}\label{variational equation-SDE}
\left\{
\begin{aligned}
d\hat{X}(r)&=\Big\{\bar{b}_x(r)\hat{X}(r)+\varepsilon_1(r)\Big\}dr\\
           &\quad+\sum\limits_{j=1}^d\Big\{\bar{\sigma}_x^j(r)\hat{X}(r)+\varepsilon_2^j(r)\Big\}dW^j(r),r\in[s,T],\\
 \hat{X}(s)&=x^1-\bar{X}^{t,x;\bar{u}}(s),
\end{aligned}
\right.
\end{equation}
and the equation for $(\hat{Y}(\cdot),\hat{Z}(\cdot))$ as
\begin{equation}\label{variational equation-BSDE}
\left\{
\begin{aligned}
-d\hat{Y}(r)&=\Big\{\bar{f}_x(r)\hat{X}(r)+\bar{f}_y(r)\hat{Y}(r)+\bar{f}_z(r)\hat{Z}(r)\\
            &\quad\ +\varepsilon_3(r)\Big\}dr-\hat{Z}(r)dW(r),\ r\in[s,T],\\
  \hat{Y}(T)&=\phi_x(\bar{X}^{t,x;\bar{u}}(T))\hat{X}(T)+\varepsilon_4(T),
\end{aligned}
\right.
\end{equation}
respectively, where
{\small\begin{equation*}
\left\{
\begin{aligned}
\varepsilon_1(r)&:=\int_0^1\big[b_x(r,\bar{X}^{t,x;\bar{u}}(r)+\theta\hat{X}(r),\bar{u}(r))\\
                &\qquad-\bar{b}_x(r)\big]\hat{X}(r)d\theta,\\
\varepsilon_2^j(r)&:=\int_0^1\big[\sigma_x^j(r,\bar{X}^{t,x;\bar{u}}(r)+\theta\hat{X}(r),\bar{u}(r))\\
                &\qquad-\bar{\sigma}_x^j(r)\big]\hat{X}(r)d\theta,\ j=1,2,\cdots,d,\\
\varepsilon_3(r)&:=\int_0^1\big[f_x(r,\bar{X}^{t,x;\bar{u}}(r)+\theta\hat{X}(r),\bar{Y}^{t,x;\bar{u}}(r)+\theta\hat{Y}(r),\\
                &\qquad\bar{Z}^{t,x;\bar{u}}(r)+\theta\hat{Z}(r),\bar{u}(r))-\bar{f}_x(r)\big]\hat{X}(r)d\theta\\
                &\quad+\int_0^1\big[f_y(r,\bar{X}^{t,x;\bar{u}}(r)+\theta\hat{X}(r),\bar{Y}^{t,x;\bar{u}}(r)+\theta\hat{Y}(r),\\
                &\qquad\bar{Z}^{t,x;\bar{u}}(r)+\theta\hat{Z}(r),\bar{u}(r))-\bar{f}_y(r)\big]\hat{Y}(r)d\theta\\
                &\quad+\int_0^1\big[f_z(r,\bar{X}^{t,x;\bar{u}}(r)+\theta\hat{X}(r),\bar{Y}^{t,x;\bar{u}}(r)+\theta\hat{Y}(r),\\
                &\qquad\bar{Z}^{t,x;\bar{u}}(r)+\theta\hat{Z}(r),\bar{u}(r))-\bar{f}_z(r)\big]\hat{Z}(r)d\theta,\\
\varepsilon_4(T)&:=\int_0^1\big[\phi_x(\bar{X}^{t,x;\bar{u}}(T)+\theta\hat{X}(T))\\
                &\quad-\phi_x(\bar{X}^{t,x;\bar{u}}(T))\big]\hat{X}(T)d\theta.
\end{aligned}
\right.
\end{equation*}}

As in pp. 258, Section 4, Chapter 5 of Yong and Zhou \cite{YZ99}, for any $k\geq1$, there exists a deterministic continuous and increasing function $\delta:[0,\infty)\rightarrow[0,\infty)$,
independent of $x^1\in\mathbf{R}^n$, with $\frac{\delta(r)}{r}\rightarrow0$ as $r\rightarrow0$, such that
\begin{equation}\label{balance estimate 1}
\left\{
\begin{aligned}
&\mathbb{E}\Big[\int_s^T|\varepsilon_1(r)|^{2k}dr\big|\mathcal{F}_s^t\Big]\leq\delta(|x^1-\bar{X}^{t,x;\bar{u}}(s)|^{2k}),\ \mathbf{P}\mbox{-}a.s.,\\
&\mathbb{E}\Big[\int_s^T|\varepsilon_2(r)|^{2k}dr\big|\mathcal{F}_s^t\Big]\leq\delta(|x^1-\bar{X}^{t,x;\bar{u}}(s)|^{2k}),\ \mathbf{P}\mbox{-}a.s.,\\
&\mathbb{E}\Big[|\varepsilon_4(T)|^{2k}\big|\mathcal{F}_s^t\Big]\leq\delta(|x^1-\bar{X}^{t,x;\bar{u}}(s)|^{2k}),\ \mathbf{P}\mbox{-}a.s.
\end{aligned}
\right.
\end{equation}
Moreover, for some $0<\alpha<1$, we have
\begin{equation}\label{balance estimate 2}
\mathbb{E}\Big[\int_s^T|\varepsilon_3(r)|^{1+\alpha}dr\big|\mathcal{F}_s^t\Big]\leq\delta(|x^1-\bar{X}^{t,x;\bar{u}}(s)|^{1+\alpha}),\ \mathbf{P}\mbox{-}a.s.
\end{equation}
In fact, denote
\begin{equation*}
\left\{
\begin{aligned}
\Delta f_x(\theta)&:=f_x(r,\bar{X}^{t,x;\bar{u}}(r)+\theta\hat{X}(r),\bar{Y}^{t,x;\bar{u}}(r)+\theta\hat{Y}(r),\\
                  &\qquad\bar{Z}^{t,x;\bar{u}}(r)+\theta\hat{Z}(r),\bar{u}(r))-\bar{f}_x(r),\\
\Delta f_y(\theta)&:=f_y(r,\bar{X}^{t,x;\bar{u}}(r)+\theta\hat{X}(r),\bar{Y}^{t,x;\bar{u}}(r)+\theta\hat{Y}(r),\\
                  &\qquad\bar{Z}^{t,x;\bar{u}}(r)+\theta\hat{Z}(r),\bar{u}(r))-\bar{f}_y(r),\\
\Delta f_z(\theta)&:=f_z(r,\bar{X}^{t,x;\bar{u}}(r)+\theta\hat{X}(r),\bar{Y}^{t,x;\bar{u}}(r)+\theta\hat{Y}(r),\\
                  &\qquad\bar{Z}^{t,x;\bar{u}}(r)+\theta\hat{Z}(r),\bar{u}(r))-\bar{f}_z(r).
\end{aligned}
\right.
\end{equation*}
Then
\begin{equation*}
\begin{aligned}
&\mathbb{E}\Big[\int_s^T|\varepsilon_3(r)|^{1+\alpha}dr\big|\mathcal{F}_s^t\Big]\\
&\leq3\mathbb{E}\Big[\int_s^T\Big|\int_0^1\Delta f_x(\theta)d\theta\Big|^{1+\alpha}|\hat{X}(r)|^{1+\alpha}dr\big|\mathcal{F}_s^t\Big]\\
&\quad+3\mathbb{E}\Big[\int_s^T\Big|\int_0^1\Delta f_y(\theta)d\theta\Big|^{1+\alpha}|\hat{Y}(r)|^{1+\alpha}dr\big|\mathcal{F}_s^t\Big]\\
&\quad+3\mathbb{E}\Big[\int_s^T\Big|\int_0^1\Delta f_z(\theta)d\theta\Big|^{1+\alpha}|\hat{Z}(r)|^{1+\alpha}dr\big|\mathcal{F}_s^t\Big]\\
&:=3I_1+3I_2+3I_3.
\end{aligned}
\end{equation*}
We consider $I_3$ only. In fact, by H\"{o}lder's inequality, for $p=\frac{2}{1-\alpha},q=\frac{2}{1+\alpha}$, we have
\begin{equation*}
\begin{aligned}
I_3&=\mathbb{E}\Big[\int_s^T\Big|\int_0^1\Delta f_z(\theta)d\theta\Big|^{1+\alpha}|\hat{Z}(r)|^{1+\alpha}dr\big|\mathcal{F}_s^t\Big]\\
   &\leq\mathbb{E}\Big[\Big(\int_s^T\Big|\int_0^1\Delta f_z(\theta)d\theta\Big|^{(1+\alpha)p}dr\Big)^{\frac{1}{p}}\\
   &\qquad\Big(\int_s^T|\hat{Z}(r)|^2dr\Big)^{\frac{1+\alpha}{2}}\big|\mathcal{F}_s^t\Big]\\
   &\leq\bigg\{\mathbb{E}\Big[\Big(\int_s^T\Big|\int_0^1\Delta f_z(\theta)d\theta\Big|^{(1+\alpha)p}dr\Big)^{\frac{2}{p}}\big|\mathcal{F}_s^t\Big]\bigg\}^{\frac{1}{2}}\\
   &\qquad\bigg\{\mathbb{E}\Big[\Big(\int_s^T|\hat{Z}(r)|^2dr\Big)^{1+\alpha}\big|\mathcal{F}_s^t\Big]\bigg\}^{\frac{1}{2}}\\
   &:=\Pi(\hat{X}(s))|x^1-\bar{X}^{t,x;\bar{u}}(s)|^{1+\alpha},
\end{aligned}
\end{equation*}
since by the second inequality of (\ref{estimate of BSDE}), we have
\begin{equation*}
\mathbb{E}\Big[\Big(\int_s^T|\hat{Z}(r)|^2dr\Big)^{1+\alpha}\big|\mathcal{F}_s^t\Big]\leq C|x^1-\bar{X}^{t,x;\bar{u}}(s)|^{2(1+\alpha)},
\end{equation*}
where
{\small\begin{equation*}
\Pi(\hat{X}(s)):=C\bigg\{\mathbb{E}\Big[\Big(\int_s^T\Big|\int_0^1\Delta f_z(\theta)d\theta\Big|^{(1+\alpha)p}dr\Big)^{\frac{2}{p}}\big|\mathcal{F}_s^t\Big]\bigg\}^{\frac{1}{2}}.
\end{equation*}}
Since from ${\bf(H3)}$ we have that $\Delta f_z(\cdot)$ is bounded and $f_z$ is continuous, then from dominate convergence theorem, we have $\Pi(\hat{X}(s))\rightarrow0$, as $x^1-\bar{X}^{t,x;\bar{u}}(s)\rightarrow0$. That is, $I_3\leq\delta(|x^1-\bar{X}^{t,x;\bar{u}}(s)|^{1+\alpha})$.

Similarly, by (\ref{estimate of SDE}) and the first inequality of (\ref{estimate of BSDE}), we can obtain the same estimates for $I_1,I_2$. Thus (\ref{balance estimate 2}) holds.

Applying It\^{o}'s formula to $\langle\hat{X}(\cdot),p(\cdot)\rangle+\hat{Y}(\cdot)q(\cdot)$, noting (\ref{adjoint equation}), (\ref{variational equation-SDE}) and (\ref{variational equation-BSDE}), we have
\begin{equation}\label{applying Ito's formula}
\begin{aligned}
&\hat{Y}(s)q(s)=-\langle\hat{X}(s),p(s)\rangle+\mathbb{E}\big[\varepsilon_4(T)q(T)\big|\mathcal{F}_s^t\big]\\
&\hspace{-2mm}-\mathbb{E}\Big[\int_s^T\langle\varepsilon_1(r),p(r)\rangle dr\big|\mathcal{F}_s^t\Big]-\mathbb{E}\Big[\int_s^T\langle\varepsilon_2(r),k(r)\rangle dr\big|\mathcal{F}_s^t\Big]\\
&\hspace{-2mm}-\mathbb{E}\Big[\int_s^T\varepsilon_3(r)q(r)dr\big|\mathcal{F}_s^t\Big],\quad\mathbf{P}\mbox{-}a.s.
\end{aligned}
\end{equation}
Noting (\ref{balance estimate 1}) and (\ref{balance estimate 2}), since $\mathbb{E}\big[\sup\limits_{s\leq r\leq T}|p(r)|^{2k}\big|\mathcal{F}_s^t\big]<\infty,\ \mathbb{E}\big[\sup\limits_{s\leq r\leq T}|q(r)|^{2k}\big|\mathcal{F}_s^t\big]<\infty,\ \mathbb{E}\Big[\int_s^T|k(r)|^{2}dr\big|\mathcal{F}_s^t\Big]<\infty$, it follows that
{\small\begin{equation*}
\begin{aligned}
&\mathbb{E}\Big[\varepsilon_4(T)q(T)\big|\mathcal{F}_s^t\Big]
\leq\Big(\mathbb{E}\Big[|\varepsilon_4(T)|^2\big|\mathcal{F}_s^t\Big]\Big)^{\frac{1}{2}}\Big(\mathbb{E}\Big[|q(T)|^2\big|\mathcal{F}_s^t\Big]\Big)^{\frac{1}{2}}\\
&\leq o(|x^1-\bar{X}^{t,x;\bar{u}}(s)|),
\end{aligned}
\end{equation*}}
{\small\begin{equation*}
\begin{aligned}
&\mathbb{E}\Big[\int_s^T\langle\varepsilon_1(r),p(r)\rangle dr\big|\mathcal{F}_s^t\Big]
\leq\mathbb{E}\Big[\sup\limits_{s\leq r\leq T}p(r)\int_s^T\varepsilon_1(r)dr\big|\mathcal{F}_s^t\Big]\\
&\leq\Big(\mathbb{E}\Big[\sup\limits_{s\leq r\leq T}|p(r)|^2\big|\mathcal{F}_s^t\Big]\Big)^{\frac{1}{2}}\Big(\mathbb{E}\Big[\Big(\int_s^T\varepsilon_1(r)dr\Big)^2\big|\mathcal{F}_s^t\Big]\Big)^{\frac{1}{2}}\\
&\leq o(|x^1-\bar{X}^{t,x;\bar{u}}(s)|),
\end{aligned}
\end{equation*}}
{\small\begin{equation*}
\begin{aligned}
&\mathbb{E}\Big[\int_s^T\langle\varepsilon_3(r),k(r)\rangle dr\big|\mathcal{F}_s^t\Big]\\
&\leq\Big(\mathbb{E}\Big[\Big(\int_s^Tk(r)dr\Big)^2\big|\mathcal{F}_s^t\Big]\Big)^{\frac{1}{2}}\Big(\mathbb{E}\Big[\Big(\int_s^T\varepsilon_3(r)dr\Big)^2\big|\mathcal{F}_s^t\Big]\Big)^{\frac{1}{2}}\\
&\leq o(|x^1-\bar{X}^{t,x;\bar{u}}(s)|),
\end{aligned}
\end{equation*}}
and
{\small\begin{equation*}
\begin{aligned}
&\mathbb{E}\Big[\int_s^T\varepsilon_3(r)q(r)dr\big|\mathcal{F}_s^t\Big]\leq\mathbb{E}\Big[\sup\limits_{s\leq r\leq T}q(r)\int_s^T\varepsilon_3(r)dr\big|\mathcal{F}_s^t\Big]\\
&\leq\Big(\mathbb{E}\Big[\sup\limits_{s\leq r\leq T}|q(r)|^q\big|\mathcal{F}_s^t\Big]\Big)^{\frac{1}{q}}\Big(\mathbb{E}\Big[\Big(\int_s^T\varepsilon_3(r)dr\Big)^{1+\alpha}\big|\mathcal{F}_s^t\Big]\Big)^{\frac{1}{1+\alpha}}\\
&\leq o(|x^1-\bar{X}^{t,x;\bar{u}}(s)|),
\end{aligned}
\end{equation*}}
where $q=\frac{1+\alpha}{\alpha}$. Thus, we have
\begin{equation}\label{y,q}
\hat{Y}(s)q(s)=-\langle\hat{X}(s),p(s)\rangle+o(|x^1-\bar{X}^{t,x;\bar{u}}(s)|),\ \mathbf{P}\mbox{-}a.s.
\end{equation}
Since $q(\cdot)$ is invertible, then
\begin{equation}\label{y}
\hat{Y}(s)=-\langle\hat{X}(s),p(s)q^{-1}(s)\rangle+o(|x^1-\bar{X}^{t,x;\bar{u}}(s)|),\ \mathbf{P}\mbox{-}a.s.
\end{equation}

Let us call a $x^1\in\mathbf{R}^n$ rational if all its coordinates are rational numbers. Since the set of all rational $x^1\in\mathbf{R}^n$ is countable, we may find a subset $\Omega_0\subseteq\Omega$ with $\mathbf{P}(\Omega_0)=1$ such that for any $\omega_0\in\Omega_0$,
\begin{equation*}
\left\{
\begin{aligned}
&V(s,\bar{X}^{t,x;\bar{u}}(s,\omega_0))=-\bar{Y}^{t,x;\bar{u}}(s,\omega_0),\\
&(\ref{estimate of SDE}), (\ref{estimate of BSDE}),(\ref{balance estimate 1}),(\ref{balance estimate 2}),(\ref{applying Ito's formula}),(\ref{y})\mbox{ are satisfied for any}\\
&\mbox{ rational }x^1,\mbox{ and }\big(\Omega,\mathcal{F},\mathbf{P}(\cdot|\mathcal{F}_s^t)(\omega_0),W(\cdot)-W(s);\\
&\ u(\cdot))|_{[s,T]}\big)\in\mathcal{U}^w[s,T].
\end{aligned}
\right.
\end{equation*}
The first equality of the above is due to the DPP (see Theorem 5.4 of Peng \cite{Peng97}). Let $\omega_0\in\Omega_0$ be fixed, then for any rational $x^1\in\mathbf{R}^n$, noting (\ref{y}), we have
\begin{equation}\label{difference}
\begin{aligned}
    &V(s,x^1)-V(s,\bar{X}^{t,x;\bar{u}}(s,\omega_0))\\
\leq&\ -Y^{s,x^1;\bar{u}}(s,\omega_0)+\bar{Y}^{t,x;\bar{u}}(s,\omega_0):=-\hat{Y}(s,\omega_0)\\
   =&\ \langle\hat{X}(s,\omega_0),p(s,\omega_0)q^{-1}(s,\omega_0)\rangle+o(|x^1-\bar{X}^{t,x;\bar{u}}(s,\omega_0)|)\\
   =&\ \langle p(s,\omega_0)q^{-1}(s,\omega_0),X^{s,x^1;\bar{u}}(s)-\bar{X}^{t,x;\bar{u}}(s,\omega_0)\rangle\\
    &\ +o(|x^1-\bar{X}^{t,x;\bar{u}}(s,\omega_0)|).
\end{aligned}
\end{equation}
Note that the term $o(|x^1-\bar{X}^{t,x;\bar{u}}(s,\omega_0)|)$ in the above depends only on the size of $|x^1-\bar{X}^{t,x;\bar{u}}(s,\omega_0)|$, and it is independent of $x^1$. Therefore, by the continuity of $V(s,\cdot)$, we see that (\ref{difference}) holds for all $x^1\in\mathbf{R}^n$, which by definition (\ref{first-order super- and sub-jets}) proves
\begin{equation*}
\big\{p(s)q^{-1}(s)\big\}\in D_x^{1,+}V(s,\bar{X}^{t,x;\bar{u}}(s)),\ \forall s\in[t,T],\ \mathbf{P}\mbox{-}a.s.
\end{equation*}
Let us now show $D_x^{1,-}V(s,\bar{X}^{t,x;\bar{u}}(s))\subset\big\{p(s)q^{-1}(s)\big\}$. Fix an $\omega\in\Omega$ such that (\ref{difference}) holds for any $x^1\in\mathbf{R}^n$. For any $\xi\in D_x^{1,-}V(s,\bar{X}^{t,x;\bar{u}}(s))$, by definition (\ref{first-order super- and sub-jets}) we have
\begin{equation*}
\begin{aligned}
0&\leq\lim\limits_{x^1\rightarrow\bar{X}^{t,x;\bar{u}}(s)}\left\{\frac{V(s,x^1)-V(s,\bar{X}^{t,x;\bar{u}}(s))}{|x^1-\bar{X}^{t,x;\bar{u}}(s)|}\right.\\
 &\qquad\qquad\qquad\quad\left.-\frac{\langle\xi,x^1-\bar{X}^{t,x;\bar{u}}(s)\rangle}{|x^1-\bar{X}^{t,x;\bar{u}}(s)|}\right\}\\
 &\leq\lim\limits_{x^1\rightarrow\bar{X}^{t,x;\bar{u}}(s)}\frac{\langle p(s)q^{-1}(s)-\xi,x^1-\bar{X}^{t,x;\bar{u}}(s)\rangle}{|x^1-\bar{X}^{t,x;\bar{u}}(s)|}.
\end{aligned}
\end{equation*}
Then, it is necessary that
\begin{equation*}
\xi=p(s)q^{-1}(s),\quad\forall s\in[t,T],\quad\mathbf{P}\mbox{-}a.s.
\end{equation*}
Thus, (\ref{connection}) holds. The proof is complete.\quad$\Box$

\vspace{1mm}

\noindent{\bf Remark 3.1}\quad Note that if $V$ is differentiable with respect to $x$, then (\ref{connection}) reduces to
\begin{equation}\label{Shi10-2}
p(s)q^{-1}(s)=V_x(s,\bar{X}^{t,x;\bar{u}}(s)),\forall s\in[t,T],\mathbf{P}\mbox{-}a.s.,
\end{equation}
which coincides with the first relation in (\ref{Shi10}) of Shi \cite{Shi10}. We point out that Theorem 3.1 is a true extension, by which we mean that it is possible to have strict set inclusions in (\ref{connection}). The following example gives such a situation.

\vspace{1mm}

\noindent{\bf Example 3.1}\quad Consider the following controlled SDE ($n=d=1$):
\begin{equation}\label{controlled SDE:example}
\left\{
\begin{aligned}
 dX^{t,x;u}(s)&=X^{t,x;u}(s)u(s)ds\\
              &\quad+X^{t,x;u}(s)dW(s),\quad s\in[t,T],\\
  X^{t,x;u}(t)&=0,
\end{aligned}
\right.
\end{equation}
with the control domain being $\mathbf{U}=[0,1]$. The cost functional is defined as
\begin{equation}\label{cost functional:example}
\begin{aligned}
&J(t,x;u(\cdot)):=-Y^{t,x;u}(s)|_{s=t},\ (t,x)\in[0,T]\times\mathbf{R}^n.
\end{aligned}
\end{equation}
with
\begin{equation}\label{controlled BSDE:example}
\left\{
\begin{aligned}
-dY^{t,x;u}(s)&=\big[X^{t,x;u}(s)-Y^{t,x;u}(s)\big]ds\\
              &\quad-Z^{t,x;u}(s)dW(s),\ s\in[t,T],\\
  Y^{t,x;u}(T)&=X^{t,x;u}(T).
\end{aligned}
\right.
\end{equation}
The corresponding generalized HJB equation reads
\begin{equation}\label{HJB equation:example}
\left\{
\begin{aligned}
 &-v_t(t,x)-\frac{1}{2}x^2v_{xx}(t,x)+x+v(t,x)\\
 &+\sup\limits_{u\in\mathbf{U}}\big\{-v_x(t,x)xu\big\}=0,\ (t,x)\in[0,T)\times\mathbf{R}^n,\\
 &v(T,x)=-x,\quad \forall x\in\mathbf{R}^n,
\end{aligned}
\right.
\end{equation}
It is not difficult to directly verify that the following function is a viscosity solution to (\ref{HJB equation:example}):
\begin{equation}
\begin{aligned}
V(t,x)=\left\{\begin{array}{lc}-x,&\mbox{if }x\leq0,\\
                               -x(T-t)-x,&\mbox{if }x>0,\end{array}\right.
\end{aligned}
\end{equation}
which obviously satisfies (\ref{value function: regularity}). Thus, by the uniqueness of the viscosity solution, $V$ coincides with the value function of our problem. Moreover, the adjoint equation writes
\begin{equation}\label{adjoint equation:example}
\left\{
\begin{aligned}
-dp(s)&=\big[\bar{u}(s)p(s)-q(s)+k(s)\big]ds-k(s)dW(s),\\
 dq(s)&=-q(s)ds,\ s\in[t,T],\\
  p(T)&=-q(T),\quad q(t)=1.
\end{aligned}
\right.
\end{equation}

Let us consider an admissible control $\bar{u}(\cdot)\equiv0$ for initial state $x=0$. The corresponding state under $\bar{u}(\cdot)$ is easily seen to be $\bar{X}^{t,x;\bar{u}}(\cdot)\equiv0$. By the stochastic verification theorem (see Theorem 9 in \cite{Zhang12}), one can check that $(\bar{X}^{t,x;\bar{u}}(\cdot),\bar{u}(\cdot))$ is really optimal. Now let us compare our main result Theorem 3.1 with the one of Shi \cite{Shi10}. In fact, by applying the results of \cite{Shi10}, especially (\ref{Shi10-2}), we obtain nothing, since $V_x(t,x)$ does not exist along the whole state $\bar{X}^{t,x;\bar{u}}(s),s\in[t,T]$. However, we have
\begin{equation}\label{connection:example}
\begin{aligned}
&D_x^{1,-}V(s,\bar{X}^{t,x;\bar{u}}(s))=\emptyset,\\
&D_x^{1,+}V(s,\bar{X}^{t,x;\bar{u}}(s))=[-(T-s)-1,-1],
\end{aligned}
\end{equation}
and the adjoint process triple is $(p(s),q(s),k(s))=(-e^{t-s},e^{t-s},0), s\in[t,T]$. Thus the relation (\ref{connection}) holds,  which shows that our Theorem 3.1 works.

\section{Concluding Remarks}

In this paper, we have established a nonsmooth version of the connection between the maximum principle and dynamic programming principle, for the stochastic recursive control problem when the control domain is convex. By employing the viscosity solution, the connection is now interpreted as a set inclusion among sub-jet  $D_x^{1,-}V(s,\bar{X}^{t,x;\bar{u}}(s))$, super-jet $D_x^{1,+}V(s,\bar{X}^{t,x;\bar{u}}(s))$ and singleton $\{p(s)q^{-1}(s)\}$. This new result has extended the classical one of Shi \cite{Shi10}, by eliminating the smoothness assumption on the value function.

This paper is the first part of our recent results on the relationship between maximum principle and dynamic programming principle under the framework of viscosity solutions, for the stochastic recursive optimal control problem. The main result in this paper (Theorem 3.1) is in local form. In the second part, we will deal with its global form, that is, the control domain is not necessarily convex. However, it looks like a difficult problem since the integrablity/regularity property of $z$ (the martingale part of the BSDE, which appears in the diffusion coefficient of the forward equation), seems to be not enough in the case when a second-order expression is necessary. In forthcoming research, we will try to overcome this difficulty by using new first- and second-order adjoint equations to deal with the global case.



\begin{thebibliography}{99}

\bibitem{Ben82}A. Bensoussan, Lectures on stochastic control, \emph{Lecture Notes in Mathematics}, vol. 972, Springer-Verlag, Berlin, 1982.

\bibitem{Bis78}J.M. Bismut, An introductory approach to duality in optimal stochastic control. \emph{SIAM Review}, {\bf 20(1)}, 62-78, 1978.

\bibitem{CZ13}J. Cvitanic, J.F. Zhang, \emph{Contract Theory in Continuous Time Models}, Springer-Verlag, Berlin, 2013.

\bibitem{DE92}D. Duffie, L.G. Epstein, Stochastic differential utility. \emph{Econometrica}, {\bf 60(2)}, 353--394, 1992.

\bibitem{EPQ97}N. El Karoui, S.G. Peng and M.C. Quenez, Backward stochastic differential equations in finance. \emph{Math. Finance}, {\bf 7(1)}, 1-71, 1997.

\bibitem{EPQ01}N. El Karoui, S.G. Peng and M.C. Quenez, A dynamic maximum principle for the optimization of recursive utilities under constraints. \emph{Ann. Appl. Proba.}, {\bf 11(3)}, 664-693, 2001.

\bibitem{FR75}W.H. Fleming, R.W. Rishel, \emph{Deterministic and Stochastic Optimal Control}, Springer-Verlag, New York, 1975.

\bibitem{N97} N. Nadirashvili, Nonuniqueness in the martingale problem and the Dirichlet problem for uniformly elliptic operators. \emph{Annali Della Scuola Normale Superiore Di Pisa Classe Di Scienze}, {\bf 24(24)}, 537-549, 1997.

\bibitem{PP90}E. Pardoux, S.G. Peng, Adapted solution of a backward stochastic differential equation. \emph{Syst. \& Control Lett.}, {\bf 14(1)}, 55-61, 1990.

\bibitem{Peng92}S.G. Peng, A generalized dynamic programming principle and Hamilton-Jacobi-Bellmen equation. \emph{Stoch. \& Stoch. Reports}, {\bf 38(2)}, 119-134, 1992.

\bibitem{Peng93}S.G. Peng, Backward stochastic differential equations and applications to the optimal control. \emph{Appl. Math. Optim.}, {\bf 27(2)}, 125-144, 1993.

\bibitem{Peng97}S.G. Peng, Backward stochastic differential equations--stochastic optimization theory and viscosity solutions of HJB equations. \emph{Topics on Stochastic Analysis}, J. Yan, S. Peng, S. Fang and L. Wu, eds., Beijing, Science Press, 85-138, 1997. (in Chinese)

\bibitem{Shi10}J.T. Shi, The relationship between maximum principle and dynamic programming principle for stochastic recursive optimal control problems and applications to finance. \emph{Proc. 29th Chinese Control Conf.}, 1535-1540, July 29-31, Beijing, China, 2010.

\bibitem{SY13}J.T. Shi, Z.Y. Yu, Relationship between maximum principle and dynamic programming for stochastic recursive optimal control problems and applications. \emph{Math. Prob. Engin.}, Vol. 2013, Article ID 285241, 12 pages.

\bibitem{WW09}G.C. Wang, Z. Wu, The maximum principle for stochastic recursive optimal control problems under partial information. \emph{IEEE Trans. Autom. Control}, {\bf 54(6)}, 1230-1242, 2009.

\bibitem{Wu13}Z. Wu, A general maximum principle for optimal control problems of forward-backward stochastic control systems. \emph{Automatica}, {\bf 49(5)}, 1473-1480, 2013.

\bibitem{Xu95}W.S. Xu, Stochastic maximum principle for optimal control problem of forward and backward system. \emph{J. Aust. Math. Soc., Ser. B}, {\bf 37(2)}, 172-185, 1995.

\bibitem{YZ99}J.M. Yong, X.Y. Zhou, \emph{Stochastic Controls: Hamiltonian Systems and HJB Equations}, Springer-Verlag, New York, 1999.

\bibitem{Zhang12}L.Q. Zhang, Stochastic verification theorem of forward-backward controlled systems for viscosity solutions. \emph{Syst. \& Control Lett.}, {\bf 61(5)}, 649-654, 2012.

\bibitem{Zhou90}X.Y. Zhou, Maximum principle, dynamic programming, and their connection in determinsitc control. \emph{J. Optim. Theory Appl.}, {\bf 65(2)}, 363-373, 1990.

\end{thebibliography}
\end{document}